\documentclass[12pt]{article}
\usepackage{amsmath,amsfonts,amsthm,amssymb}
\newcommand{\F}{\mathbb F_q}
\newcommand{\K}{\overline{K}_c}

\begin{document}
\newtheorem{prop}{Proposition}
\pagestyle{plain}
\title{Dwork-Carlitz Exponential and Overconvergence for Additive 
Functions in Positive Characteristic}
\author{Anatoly N. Kochubei\footnote{Partially supported by
DFG under Grant 436 UKR 113/87/01, and by the Ukrainian Foundation 
for Fundamental Research, Grant 10.01/004.}\\ 
\footnotesize Institute of Mathematics,\\ 
\footnotesize National Academy of Sciences of Ukraine,\\ 
\footnotesize Tereshchenkivska 3, Kiev, 01601 Ukraine
\\ \footnotesize E-mail: \ kochubei@i.com.ua}
\date{}
\maketitle
\vspace{3cm}
\begin{abstract}
We study overconvergence phenomena for $\F$-linear functions on a 
function field over a finite field $\F$. In particular, an analog 
of the Dwork exponential is introduced.\\

\end{abstract}

\newpage

{\bf 1.} One of the central subjects of contemporary $p$-adic 
analysis is that of overconvergence. In contrast to analysis over 
$\mathbb R$ and $\mathbb C$, the power series for principal 
special functions over $\mathbb Q_p$ or $\mathbb C_p$ converge 
only on finite disks or annuli. For example, the exponential 
series $\exp (t)=\sum\limits_{n=0}^\infty \frac{t^n}{n!}$, $t\in 
\mathbb C_p$, converges if and only if $|t|_p<p^{-1/(p-1)}$ (see 
\cite{DGS} or \cite{R}).

At the same time, for many special functions there exist some 
expressions combining their values in various points (usually 
connected by the Frobenius power $t\mapsto t^p$), for which the 
corresponding power series converge on wider regions. The simplest 
example is the Dwork exponential
\begin{equation}
\theta (t)=\exp (\pi (t-t^p))
\end{equation}
where $\pi$ is a root of the equation $z^{p-1}+p=0$. The power 
series for $\theta (t)$, in the variable $t$, converges for 
$|t|_p<p^{\frac{p-1}{p^2}}$ ($>1$), though the formula (1) is not 
valid outside the unit disk. The special value $\theta (1)$ is a 
primitive $p$-th root of unity.

Other examples involve the exponential function of $q$-analysis 
\cite{AV}, some hypergeometric functions \cite{D}, polylogarithms 
\cite{Col}, and many others. The overconvergent functions usually 
satisfy equations possessing special algebraic properties called 
the Frobenius structures (see \cite{AV,RC}).

In this paper we consider the overconvergence phenomena in the 
case of a local field of a positive characteristic, that is (up to 
an isomorphism) the field $K$ of formal Laurent series
\begin{equation}
z=\sum\limits_{i=n}^\infty \zeta_ix^i,\quad n\in \mathbb 
Z,\ \zeta_i\in \F ,
\end{equation}
with coefficients from a finite field $\F$. If $z\in K$ is an 
element (2) with $\zeta_n\ne 0$, its non-Archimedean absolute 
value $|z|$ is given by $|z|=q^{-n}$. Note that the above 
construction of $K$ (in contrast to the completion of $\F (x)$ 
with respect to the $\infty$-valuation often used in the function 
field arithmetic) leads to structures (operators, orthonormal 
bases etc) resembling to some extent the $p$-adic case; see, for example, 
\cite{K98,Con}.

We will consider only $\F$-linear functions and power series of 
the form $\sum\limits_{k=0}^\infty c_kt^{q^k}$, $c_k\in \K$ ($\K$ 
is the completion of an algebraic closure $\overline{K}$ of $K$). 
This class contains many important functions, including the 
Carlitz exponential and logarithm, analogs of the Bessel and 
hypergeometric functions, polylogarithms etc.

In particular, using the Carlitz exponential $e_C$ (see below), we 
construct an analog of the Dwork exponential and prove its 
overconvergence, consider the overconvergence problems for other 
special functions mentioned above. These problems are much 
simpler than those in the characteristic zero case. The reason is that 
the above functions satisfy differential equations with the 
Carlitz derivatives (see \cite{K00, K03,KPoly,KHyp,Th}); the 
difference structure of the latter leads immediately to 
overconvergence properties of some linear combinations of 
solutions.

{\bf 2.} The Carlitz exponential (see \cite{G,Th}) is the function
\begin{equation}
e_C(t)=\sum\limits_{n=0}^\infty \frac{t^{q^n}}{D_n}
\end{equation}
where $D_n$ is the Carlitz factorial
\begin{equation}
D_n=[n][n-1]^q\ldots [1]^{q^{n-1}},\quad [n]=x^{q^n}-x\ (n\ge 1),\ 
D_0=1.
\end{equation}
Since $|[n]|=q^{-1}$ for any $n\ge 1$, it follows from (4) that
$$
|D_n|=q^{-\frac{q^n-1}{q-1}},
$$
so that the series in (3) converges for $|t|<q^{-\frac{1}{q-1}}$.

Let $\sigma$ be an arbitrary solution of the equation 
$z^{q-1}=-x$. Let us consider the function
\begin{equation}
E(t)=e_C(\sigma (t-t^q)),
\end{equation}
defined initially for $|t|<q^{-\frac{1}{q-1}}$ (we denote by 
$|\cdot |$ also the extension of the absolute value from $K$ onto 
$\K$). Note that, in spite of a formal resemblance, the formulas 
for the Dwork exponential (1) and ``the Dwork-Carlitz 
exponential'' (5) have a quite different meaning -- the function 
$\theta (t)$ is a multiplicative combination of values of the 
classical exponential, while $E(t)$ is an additive combination of 
values of the Carlitz exponential. This difference from the 
classical overconvergence theory appears also in some other 
examples given below.

From (3) and (5), after a simple transformation we find that
\begin{equation}
E(t)=\sigma t+\sum\limits_{n=1}^\infty \left( 
\frac{\sigma^{q^n}}{D_n}-\frac{\sigma^{q^{n-1}}}{D_{n-1}}\right) 
t^{q^n}.
\end{equation}
In order to investigate the convergence of the series (6), we have 
to study the structure of elements $D_n$.

\medskip
\begin{prop}
For any $n\ge 1$,
\begin{equation}
\left| D_n-(-1)^nx^{1+q+\cdots +q^{n-1}}\right| \le 
q^{-\frac{q^n-1}{q-1}-(q-1)q^{n-1}}.
\end{equation}
\end{prop}

\medskip
{\it Proof}. We will prove that
\begin{equation}
\left| D_n-(-1)^nx^{1+q+\cdots +q^{n-1}}\right| \le 
q^{-l_n}
\end{equation}
where the sequence $\{ l_n\}$ is determined by the requrrence
\begin{equation}
l_n=ql_{n-1}+1,\quad l_1=q.
\end{equation}

Indeed, if $n=1$, then $D_1=x^q-x$, so that $|D_1+x|=q^{-q}$. 
Suppose that we have proved (8) for some value of $n$. We have
$$
\left| D_n^q-(-1)^nx^{q+q^2\cdots +q^n}\right| \le 
q^{-ql_n},
$$
whence
$$
\left| [n+1]D_n^q-(-1)^n[n+1]x^{q+q^2\cdots +q^n}\right| \le 
q^{-(ql_n+1)}.
$$
Since $D_{n+1}=[n+1]D_n^q$, we find that
\begin{equation}
\left| \left( D_{n+1}-(-1)^{n+1}x^{1+q+\cdots +q^n}\right)
-(-1)^nx^{q+q^2\cdots +q^{n+1}}\right| \le 
q^{-l_{n+1}}.
\end{equation}

It is easy to check that
\begin{equation}
l_n=\frac{q^n-1}{q-1}+(q-1)q^{n-1}
\end{equation}
satisfies (9); the expression (11) can also be deduced from a 
general formula for a solution of a difference equation; see 
\cite{Gel}.

On the other hand,
\begin{equation}
\left| (-1)^nx^{q+q^2\cdots +q^{n+1}}\right| 
=q^{-\frac{q^{n+2}-q}{q-1}},
\end{equation}
and, by a simple computation,
\begin{equation}
\frac{q^{n+2}-q}{q-1}-l_{n+1}=q^n-1>0,\quad n\ge 1.
\end{equation}

It follows from (10), (12), (13), and the ultra-metric property of 
the absolute value, that
$$
\left| D_{n+1}-(-1)^{n+1}x^{1+q+\cdots +q^n}\right| \le 
q^{-l_{n+1}},
$$
which proves the inequalities (8) and (7) for any $n$. \qquad 
$\blacksquare$

\medskip
Now we can prove the overconvergence of $E(t)$.

\medskip
\begin{prop}
The series in (6) converges for $|t|<\rho$, where $\rho 
=q^{\frac{q-1}{q^2}}>1$. In particular, $E(1)=\lim\limits_{n\to 
\infty }\frac{\sigma^{q^n}}{D_n}$ is defined, and $E(1)=\sigma$.
\end{prop}

\medskip
{\it Proof}. Let us write
$$
\frac{\sigma^{q^n}}{D_n}-\frac{\sigma^{q^{n-1}}}{D_{n-1}}
=\frac{\sigma^{q^{n-1}}}{D_{n-1}}\left( 
\frac{\sigma^{q^n-q^{n-1}}D_{n-1}}{D_n}-1\right) .
$$
We have $\sigma^{q^n-q^{n-1}}=-x^{q^{n-1}}$, so that
\begin{multline*}
\frac{\sigma^{q^n}}{D_n}-\frac{\sigma^{q^{n-1}}}{D_{n-1}}
=-\frac{\sigma^{q^{n-1}}}{D_{n-1}D_n}\left( x^{q^{n-1}}D_{n-1}
+D_n\right) \\
=-\frac{\sigma^{q^{n-1}}}{D_{n-1}D_n}\left\{ x^{q^{n-1}}\left( 
D_{n-1}-(-1)^{n-1}x^{1+\cdots +q^{n-2}}\right) 
+\left( D_n-(-1)^n x^{1+\cdots +q^{n-1}}\right) \right\}.
\end{multline*}

If $n\ge 2$, then by Proposition 1,
$$
\left|  x^{q^{n-1}}\left( 
D_{n-1}-(-1)^{n-1}x^{1+\cdots +q^{n-2}}\right) \right| \le
q^{-\left( q^{n-1}+\frac{q^{n-1}-1}{q-1}+(q-1)q^{n-2}\right) },
$$
$$
\left| D_n-(-1)^n x^{1+\cdots +q^{n-1}}\right| \le 
q^{-\frac{q^n-1}{q-1}-(q-1)q^{n-1}}.
$$
Comparing the right-hand sides we check that the first of them 
is bigger; therefore 
$$
\left| \frac{\sigma^{q^n}}{D_n}-\frac{\sigma^{q^{n-1}}}{D_{n-1}}
\right| \le q^{-\frac{q^{n-1}}{q-1}}\cdot 
q^{-\frac{q^{n-1}-1}{q-1}} \cdot q^{-\frac{q^n-1}{q-1}}\cdot
q^{-\left( q^{n-1}+\frac{q^{n-1}-1}{q-1}+(q-1)q^{n-2}\right)} ,
$$
so that
\begin{equation}
\left| \frac{\sigma^{q^n}}{D_n}-\frac{\sigma^{q^{n-1}}}{D_{n-1}}
\right| \le q^{-\frac{q^{n-2}(q-1)^2+1}{q-1}},\quad n\ge 2.
\end{equation}

For $n=1$, we get
$$
\left| \frac{\sigma^q}{D_1}-\sigma \right| =|\sigma |\left| 
\frac{-x}{[1]}-1\right| =|\sigma |\left| \frac{x^q}{[1]}\right| ,
$$
whence
\begin{equation}
\left| \frac{\sigma^q}{D_1}-\sigma \right| \le 
q^{-\frac{1}{q-1}-(q-1)}.
\end{equation}

It follows from (14) that the series in (6) converges for 
$|t|<\rho$. For $t=1$, we obtain that
\begin{equation}
E(1)=\sigma +\sum\limits_{n=1}^\infty \left( 
\frac{\sigma^{q^n}}{D_n}-\frac{\sigma^{q^{n-1}}}{D_{n-1}}\right) 
=\lim\limits_{n\to \infty }\frac{\sigma^{q^n}}{D_n}.
\end{equation}
Note that
$$
\left| \frac{\sigma^{q^n}}{D_n}\right| =q^{-\frac{1}{q-1}}
$$
for all values of $n$. Now
$$
E(1)^q=\lim\limits_{n\to \infty }\frac{\sigma^{q^{n+1}}}{D_n^q}
=\lim\limits_{n\to \infty }[n+1]\frac{\sigma^{q^{n+1}}}{D_{n+1}}
=\lim\limits_{n\to \infty }[n]\frac{\sigma^{q^n}}{D_n}=-xE(1)
$$
because
$$
\left| \frac{x^{q^n}\sigma^{q^n}}{D_n}\right| 
=q^{-\frac{1}{q-1}-q^n}\longrightarrow 0,
$$
as $n\to \infty$.

By (16), $E(1)\ne 0$, so that $E(1)^{q-1}=-x$, thus $E(1)$ 
satisfies the same equation as $\sigma$. All the solutions of this 
equation are obtained by multiplying $\sigma$ by non-zero elements 
$\xi \in \F$. Therefore $E(1)=\sigma \xi$, $\xi \in \F$, $\xi \ne 
0$. If $\xi \ne 1$, then
\begin{equation}
|E(1)-\sigma |=|(1-\xi )\sigma |=|\sigma |=q^{-\frac{1}{q-1}}.
\end{equation}
On the other hand, by (16),
$$
|E(1)-\sigma |\le \sup\limits_{n\ge 1}\left| 
\frac{\sigma^{q^n}}{D_n}-\frac{\sigma^{q^{n-1}}}{D_{n-1}}
\right| ,
$$
and we see that (17) contradicts (14) and (15). $\qquad 
\blacksquare$

\medskip
It is interesting that the special value $\sigma =E(1)$, just as 
the special value of the Dwork exponential in the characteristic 0 
case, generates a cyclotomic extension of the function field 
(related in this case to the Carlitz module); see \cite{Ros}.

\medskip
{\bf 3.} As it has been mentioned, many important $\F$-linear 
functions defined on subsets of $K$ satisfy equations involving 
the Carlitz derivative $d=\sqrt[q]{}\circ \Delta$ where
$$
\Delta u(t)=u(xt)-xu(t).
$$

The Carlitz exponential $e_C$ satisfies the simplest equation 
$de_C=e_C$, so that
\begin{equation}
e_C(t)^q+xe_C(t)=e_C(xt).
\end{equation}
The right-hand side of (18) obviously converges on a wider disk 
than $e_C$ itself (note that, in contrast to the $p$-adic case, 
$E(t)$ does not satisfy a homogeneous equation with the Carlitz 
derivative).

Similarly, the Bessel-Carlitz function $J_n(t)$, introduced in 
\cite{C}, satisfies the identity $\Delta J_n=J_{n-1}^q$, so that
\begin{equation}
J_{n-1}^q(t)+xJ_n(t)=J_n(xt),
\end{equation}
and we have an overconvergence for the right-hand side of (19). In 
this sense equations with the Carlitz derivatives may be seen 
themselves as analogs of the Frobenius structures of $p$-adic 
analysis.

The next two examples (of an essentially similar nature) are just 
a little more complicated.

\medskip
{\bf 4.} Polylogarithms on $K$, in the sense of \cite{KPoly}, are 
defined as follows. First the function $l_1(t)$, an analog of the 
function $-\log (1-t)$, is introduced as a solution of the 
equation $(1-\tau )du(t)=t$, where $\tau u=u^q$. This 
equation is, of course, an analog of the classical equation 
$(1-t)u'(t)=1$ (note that the 
function $f(t)=t$ is the unit element in the composition rings of 
$\F$-linear polynomials or holomorphic functions). Then the 
polylogarithms $l_n(t)$ are defined reqursively by the equations 
$\Delta l_n=l_{n-1}$, $n\ge 2$ (classically, 
$tl_n'(t)=l_{n-1}(t)$). These definitions lead \cite{KPoly} to the 
explicit expressions
\begin{equation}
l_n(t)=\sum\limits_{j=1}^\infty \frac{t^{q^j}}{[j]^n},\quad 
n=1,2,\ldots .
\end{equation}
The series in (20) converges for $|t|<1$. In \cite{KPoly} we 
constructed their continuous extensions to the disk $\{ t\in K:\ 
|t|\le 1\}$. Here we give the following overconvergence result 
resembling Coleman's theorem \cite{Col} about classical 
polylogarithms.

\medskip
\begin{prop}
The power series for the function $L_n(t)=l_n(t)-l_n(t^q)$ 
converges for $|t|<q^{1/q}$.
\end{prop}

\medskip
{\it Proof}. By a simple transformation, we get
\begin{equation}
L_n(t)=\frac{t^q}{[1]^n}+\sum\limits_{j=2}^\infty \left( 
\frac{1}{[j]^n}-\frac{1}{[j-1]^n}\right) t^{q^j}.
\end{equation}

We have,
$$
\frac{1}{[j]^n}-\frac{1}{[j-1]^n}=
\frac{([j-1]-[j])\left( [j]^{n-1}+[j-1][j]^{n-2}+\cdots 
+[j-1]^{n-1}\right) }{[j]^n[j-1]^n}.
$$
For any $j\ge 2$, $|[j]|=q^{-1}$, $|[j-1]-[j]|=\left| 
x^{q^{j-1}}-x^{q^j}\right|=q^{-q^{j-1}}$, so that
$$
\left| \frac{1}{[j]^n}-\frac{1}{[j-1]^n}\right| \le
q^{-q^{j-1}+n-1},
$$
and the convergence radius of the series (21) equals $q^{1/q}$. 
$\qquad \blacksquare$

\medskip
{\bf 5.} Let us consider the hypergeometric function \cite{KHyp}
\begin{equation}
F(a,b;c;t)=\sum\limits_{n=0}^\infty \frac{\langle a\rangle_n
\langle b\rangle_n}{\langle c\rangle_nD_n}t^{q^n},
\end{equation}
where $a,b,c\in \K$, $c\notin \{ [0],[1],\ldots ,[\infty ]\}$, 
$[\infty ]=-x$, and the Pochhammer-type symbols are defined as 
$\langle a\rangle_0=1$,
$$
\langle a\rangle_m=([0]-a)^{q^n}([1]-a)^{q^{n-1}}\cdots 
([n-1]-a)^q,\quad n\ge 1.
$$
If all the parameters have the form $[-\alpha ],\alpha \in \mathbb 
Z$, then the function (22) coincides, up to a change of variable, 
with the hypergeometric function introduced by Thakur \cite{Th}.

Denote $T_1(a)=(a-[1])^{1/q}$, $a\in \K$. The transformation $T_1$ 
is an analog of the unit shift of integers: if $a=[-\alpha ],\alpha \in \mathbb 
Z$, then $T_1([-\alpha ])=[-\alpha -1]$. The identity
\begin{equation}
\langle a\rangle_n=-a^{q^n}\langle T_1(a)\rangle_{n-1}^q,\quad 
n\ge 1,
\end{equation}
holds for any $a\in \K$ (see \cite{KHyp}).

If $|a|=|b|=|c|=1$, then $|T_1(a)|=|T_1(b)|=|T_1(c)|=1$, and the 
disk of convergence of the series (22) is the same as the one for 
the Carlitz exponential, that is $\left\{ t\in \K:\ 
|t|<q^{-\frac{1}{q-1}}\right\}$.

\medskip
\begin{prop}
The identity
\begin{equation}
\tau F(T_1(a),T_1(b);T_1(c);\frac{ab}c t)-xF(a,b;c;t)=-F(a,b;c;xt)
\end{equation}
holds for any values of the variable and parameters, such that all 
the terms of (24) make sense. In particular, if $|a|=|b|=|c|=1$, then the
right-hand side of (24) is overconvergent, that is the series for the
right-hand side converges for $|t|<q^{1-\frac{1}{q-1}}$ ($>1$).
\end{prop}

\medskip
{\it Proof}. Changing the index of summation we find that
$$
\tau F(T_1(a),T_1(b);T_1(c);z)=\sum\limits_{n=0}^\infty \frac{\langle 
T_1(a)\rangle^q_{n-1}\langle T_1(b)\rangle^q_{n-1}}{\langle 
T_1(c)\rangle^q_{n-1}D^q_{n-1}}t^{q^n} 
$$
for any $z$ from the convergence disk. Using the identity (23) and 
the fact that $D_n=[n]D_{n-1}^q$ we get
$$
\tau F(T_1(a),T_1(b);T_1(c);z)=-\sum\limits_{n=0}^\infty \frac{\langle a\rangle_n
\langle b\rangle_n[n]}{\langle c\rangle_nD_n}\left( \frac{ab}{c}\right)^{-q^n}z^{q^n}
$$
(note that $[0]=0$), which implies (24). $\qquad \blacksquare$

\end{document}